\documentclass[12pt]{amsart}

\newtheorem{theorem}{Theorem}
\newtheorem{proposition}{Proposition}
\newtheorem{lemma}{Lemma}
\newtheorem{corollary}{Corollary}

\begin{document}

\title[Free orthogonal quantum groups]{The representation theory of free orthogonal quantum groups}

\author{Teodor Banica}
\address{Department of Mathematics, Paris 7 University, 75005 Paris, France}

\begin{abstract}
We find, for each $n\geq2$, the class of $n\times n$ compact quantum groups whose representation theory is similar to that of $SU(2)$: this is the class of ``free analogues of $O(n)$'' constructed by Van Daele and Wang.
\end{abstract}

\maketitle

Let us first recall the definition of the compact matrix quantum groups \cite{wo1}, \cite{wo3}: these are the pairs $(G,u)$ formed by a unital $C^*$-algebra $G$ and a matrix $u\in M_n(G)$ such that:
\begin{enumerate}
\item The coefficients of $u$ generate a dense $*$-subalgebra $G_s\subset G$.

\item There exists a $C^*$-morphism $\delta:G\to G\otimes G$ such that $(Id\otimes\delta)u=u_{12}u_{13}$.

\item The matrices $u$ and $\bar{u}$ are both invertible.
\end{enumerate}

We call representation of $(G,u)$ any invertible matrix $r\in M_k(G_s)$ such that $(Id\otimes\delta)=r_{12}r_{13}$. In \cite{wo1} Woronowicz shows that any compact matrix quantum group has a Haar measure, and develops a Peter-Weyl type theory for its representations. We will freely use the notations and results from \cite{wo1}.

\medskip

Wang \cite{wan}, then Van Daele and Wang \cite{vwa} have recently constructed compact matrix quantum groups having universality properties similar to those of $U(n)$ and $O(n)$:

\medskip

{\bf The unitary case.} Let $(G,u)$ be a compact matrix quantum group. Then any representation of $G$ is equivalent to a unitary representation. We can therefore suppose (up to similarity) that $u$ and $F\bar{u}F^{-1}$ are unitaries, for a certain scalar matrix $F$.

We can define, for any $n\in\mathbb N$ and any $F\in GL(n)$, the universal $C^*$-algebra $A_u(F)$ generated by variables $(u_{ij})_{1\leq i,j\leq n}$, with the relations making unitaries the matrices $u$ and $F\bar{u}F^{-1}$. It is easy to see that $(A_u(F),u)$ is a compact matrix quantum group.

\medskip

{\bf The orthogonal case.} Let $(G,u)$, with $u$ being unitary, and assume that we have $u\sim\bar{u}$ as representations, so that there exists a scalar matrix $F$ such that $u=F\bar{u}F^{-1}$. We have then $\bar{u}=\bar{F}u\bar{F}^{-1}$, and so $u=(F\bar{F})u(F\bar{F})^{-1}$. Thus if $u$ is irreducible, then $F\bar{F}=c\in\mathbb R$ (because $F\bar{F}=c\in\mathbb C\implies \bar{F}F=\bar{c}\implies c\in\mathbb R$).

We can define, for any $n\in\mathbb N$ and any $F\in GL(n)$ satisfying $F\bar{F}=c\in\mathbb R$, the universal $C^*$-algebra $A_o(F)$ generated by variables $(u_{ij})_{1\leq i,j\leq n}$, with the relations $u=F\bar{u}F^{-1}=$ unitary. It is easy to see that $(A_o(F),u)$ is a compact matrix quantum group.

\medskip

{\bf Remark.} We have $A_o(^0_{-1}{\ }^1_0)=C(SU(2))$. In fact, it is easy to see that $S_\mu U(2)=A_o(^0_{-\mu^{-1}}{\ }^1_0)$, and that any $A_o(F)$ with $F\in GL(2)$ is similar to a certain $S_\mu U(2)$ (see \cite{vwa}).

\medskip

The above link between $SU(2)$ and the algebras $A_o(F)$ can be extended at $n\geq3$:

\begin{theorem}
Let $n\in\mathbb N$, and let $F\in GL(n)$ such that $F\bar{F}=c\in\mathbb R$. Then the irreducible representations of $A_o(F)$ are self-adjoint, and can be indexed by $\mathbb N$, with $r_0=1$, $r_1=u$ and
$$r_kr_s=r_{|k-s|}+r_{|k-s|+2}+\ldots+r_{k+s-2}+r_{k+s}$$
(i.e. the same formulae as for the representations of $SU(2)$).
\end{theorem}

In order to prove this result, we begin with some considerations regarding the concrete complete monoidal $W^*$-category (see \cite{wo2}) $X(F)$ of representations of $A_o(F)$. In terms of monoidal categories, the relation $u\sim\bar{u}$ tells us that:

-- $X(F)$ is the completion of the subcategory $Y(F)$ having as objects $1,u,u^2,u^3,\ldots$

-- $Y(F)$ contains a certain morphism, intertwining $1$ and $u^2$.

These two conditions allow one to fully reconstruct $Y(F)$ (see also the corresponding construction for $S_\mu U(N)$ from \cite{wo2}):

\begin{proposition}
Let $H=\mathbb C^n$, with standard basis $\{e_i\}$. For $r,s\in\mathbb N$ we define the sets $Mor(r,s)\subset B(H^{\otimes r},H^{\otimes s})$ of linear combinations of (composable) products of maps of type $Id_{H^{\otimes k}}$ or $Id_{H^{\otimes k}}\otimes E\otimes Id_{H^{\otimes p}}$ or $Id_{H^{\otimes k}}\otimes E^*\otimes Id_{H^{\otimes p}}$, where $E\in Mor(0,2)$ is the linear map $1\to\sum F_{ji}e_i\otimes e_j$. Then, each $Mor(r,s)$ equals $Mor(u^r,u^s)\subset B(H^{\otimes r},H^{\otimes s})$.
\end{proposition}

\begin{proof}
$Z(F)=\{\mathbb N,+,\{H^{\otimes r}\}_{r\in\mathbb N},\{Mor(r,s)\}_{r,s\in\mathbb N}\}$ is clearly a concrete monoidal $W^*$-category, generated by 1. If we denote by $k:H\to H$ the antilinear involution given by $\lambda e_i\to\bar{\lambda}e_i$ and $j=k\bar{F}$, then $t_j,t_{j^{-1}}\in Mor(0,2)$, and so $1=\bar{1}$ inside $Z(F)$ (see \cite{wo2}, page 39). By the duality theorem (Theorem 1.3 in \cite{wo2}) the universal $Z(F)$-admissible pair is a compact matrix quantum group defined by the same universal property as $A_o(F)$, and so is $A_o(F)$ itself. Thus $Y(F)=Z(F)$, and so $Mor(u^r,u^s)=Mor(r,s)$, for any $r,s\in\mathbb N$.
\end{proof}

In order to compute the above monoidal $W^*$-category, we will need:

\begin{lemma}
$(E^*\otimes Id_H)(Id_H\otimes E)=cId_H$.
\end{lemma}

\begin{proof}
This is clear from the definition of $E$, and from $F\bar{F}=c$.
\end{proof}

We study now the spaces $Mor(k,k)$. For $s=1,\ldots,k-1$ let us define:
$$f_s=||E(1)||^{-2}Id_{H^{\otimes s-1}}\otimes EE^*\otimes Id_{H^{\otimes k-s-1}}$$

An elementary computation based on Lemma 1 above shows that:
\begin{enumerate}
\item $f_s^2=f_s^*=f_s$, $\forall 1\leq s\leq k-1$.

\item $f_sf_t=f_tf_s$, $\forall 1\leq s,t\leq k-1$ with $|s-t|\geq2$.

\item $\beta f_sf_tf_s=f_s$, $\forall 1\leq s,t\leq k-1$ with $|s-t|=1$, with $\beta=c^{-2}||E(1)||^4$.
\end{enumerate}

We recall that the Temperley-Lieb algebra $A_{\beta,k}$ is defined with generators $1,f_1,\ldots,f_{k-1}$ and the above relations (see \cite{jon}). In order to make the link with $A_{\beta,k}$, we will need:

\begin{proposition}
The elements $1,f_1,\ldots,f_{k-1}$ generate $Mor(k,k)$, as a $\mathbb C$-algebra.
\end{proposition}

\begin{proof}
Let $I(p)=Id_{H^{\otimes p}}$ and $V(p,q)=I(p)\otimes E\otimes I(q)$. By using Lemma 1, we can see that any morphism of $Y(F)$ appears as a linear combination of maps of type $I(.)$ or of type $V(.,.)\circ\ldots\circ V(.,.)\circ V^*(.,.)\circ\ldots\circ V^*(.,.)$. In particular, the elements of $Mor(k,k)$ are linear combinations of $I(k)$ and of maps of the following form:
$$(*)\qquad V(p_m,q_m)\circ\ldots\circ V(p_1,q_1)\circ V^*(r_1,s_1)\circ\ldots\circ V^*(r_m,s_m)$$

Let $U$ be the set of morphisms of $Y(F)$ which are linear combinations of maps of the form $I(m)$ or of the form $U(m,q,p)=I(m)\otimes E\otimes I(q)\otimes E^*\otimes I(p)$ or of the form $U(m,q,p)^*$. We show by recurrence that any map of the form $(*)$ belongs to $U$.

Indeed, let us suppose that this is true at $m\in\mathbb N$, and pick an arbitrary map of the form $A=V(p_{m+1},q_{m+1})\circ\ldots\circ V(p_1,q_1)\circ V^*(r_1,s_1)\circ\ldots\circ V^*(r_{m+1},s_{m+1})$. By recurrence, we have $A=V(p_{m+1},q_{m+1})\circ T\circ V^*(r_{m+1},s_{m+1})$, for a certain $T\in U$.

It is clear that any product of the form $V(.,.)\circ U(.,.,.)^\sigma$, with $\sigma\in\{1,*\}$, can be written in the form $U(.,.,.)^\sigma\circ V(.,.)$. By using a recurrence, we conclude that any product of the form $V(.,.)\circ T$ with $T\in U$ can be written as a sum $\sum T_i\circ V(a_i,b_i)$ with $T_i\in U$. Thus $A$ is of the form $\sum T_i\circ V(.,.)\circ V^*(.,.)$, with $T_i\in U$. But, each product of the type $V(.,.)\circ V^*(.,.)$ being in $U$, we conclude that we have $A\in U$.

Summarizing, we have proved that $Mor(k,k)$ is contained in $U$, and so is the set of linear combinations of products of maps of type $I(k)$ or $U(m,p,q)$ or $U(m,q,p)^*$, with $m+q+p=k-2$. Now since $f_1,\ldots,f_{k-1}$ are self-adjoint, it remains to prove that $U(m,p,q)$ with $m+q+p=k-2$ belongs to the algebra generated by $f_1,\ldots,f_{k-1}$. But this is clear from the formula $c^q||E(1)||^{-2q-2}U(m,p,q)=f_{m+1}f_{m+2}\ldots f_{m+q}f_{m+q+1}$, which can be shown by recurrence on $q$, by using Lemma 1.
\end{proof}

As a consequence of the above result, we obtain:

\begin{corollary}
$Mor(k,k)$ is a quotient of $A_{\beta,k}$.
\end{corollary}

{\bf Remark.} By counting the reduced words in $A_{\beta,k}$ we have $\dim(A_{\beta,k})\leq C_k:=\frac{1}{k+1}\binom{2k}{k}$, with these latter numbers being the Catalan numbers (see \cite{jon}, Aside 4.1.4). Thus, we have $\dim(Mor(u^k,u^k))=\dim(Mor(k,k))\leq\dim(A_{\beta,k})\leq C_k$, for any $k$. 

In the particular case of $A_o(^0_{-1}{\ }^1_0)=C(SU(2))$ we have equality everywhere. This is well-known, and can be proved as well as follows. Let $v$ be the fundamental representation of $SU(2)$, with character denoted $f$, and let $h$ be the Haar integration over $SU(2)$. The, the last formula from the first appendix in \cite{wo1} gives $h((f/2)^{2k})=\frac{2}{\pi}\int_{-1}^1x^{2k}\sqrt{1-x^2}dx$. In addition, $h((f/2)^{2k+1})=0$, and so $f/2\in(C(SU(2)),h)$ is a semicircular variable in the sense of Voiculescu \cite{vdn}. Thus $dim(Mor(v^k,v^k))=4^kh((f/2)^{2k})=4^k\gamma_{0,1}(X^{2k})$, where $\gamma_{0,1}$ is the semicircle law, whose moments can be computed by using 3.3 and 3.4 in \cite{vdn} and the residue formula, as follows: $\gamma_{0,1}(X^{2k})=\frac{1}{2k+1}\cdot\frac{1}{2\pi i}\int_T(z^{-1}+z/4)^{2k+1}=4^{-k}C_k$.

\medskip

Now back to our quantum group setting, we obtain from this:

\begin{corollary}
Let $u,v$ be the fundamental representations of $A_o(F),SU(2)$, respectively. Then $\dim(Mor(u^k,u^k))\leq\dim(Mor(v^k,v^k))$.
\end{corollary}

With these results in hand, we can now prove Theorem 1:

\begin{proof}(of Theorem 1) Let $\{\chi_k\}_{k\in\mathbb N}$ be the characters of the irreducible representations of $SU(2)$. The linear space $A\subset C(SU(2))$ spanned by these characters is then a $\mathbb C$-algebra, which is isomorphic to $\mathbb C[X]$, via $X\to\chi_1$. By recurrence on $k$, we can find integers $a(k,s)\in\mathbb N$ such that $a(k,k)=1$ and $\chi_1^k=\sum_{s=0}^ka(k,s)\chi_s$.

Since $A$ is a polynomial algebra on $\chi_1$, we can define a morphism $\Psi:A\to A_o(F)$ by $\chi_1\to f_1$, where $f_1$ is the character of the fundamental representation of $A_o(F)$. The elements $f_k=\Psi(\chi_k)\in A_o(F)$ verify then $f_kf_s=f_{|k-s|}+f_{|k-s|+2}+\ldots+f_{k+s}$.

We show now by recurrence on $k$ that each $f_k$ is the character of an irreducible representation $r_k$ of $A_o(F)$, non-equivalent to $r_0,\ldots,r_{k-1}$. At $k=0,1$ this is clear. Assume now that the result holds at $k-1$. We have $f_{k-2}f_1=f_{k-3}+f_{k-1}$, and so $r_{k-2}r_1=r_{k-3}+r_{k-1}$, which gives $r_{k-1}\subset r_{k-2}r_1$. Now since $r_{k-2}$ is by recurrence irreducible, by Frobenius reciprocity we have $r_{k-2}\subset r_{k-1}r_1$, so there exists a representation $r_k$ such that $r_{k-1}r_1=r_{k-2}+r_k$. Since $f_{k-1}f_1=f_{k-2}+f_k$, the character of $r_k$ is $f_k$.

Now since $\Psi$ is a morphism, we have $f_1^k=\sum_{s=0}^ka(k,s)f_s$ and so $\dim(Mor(u^k,u^k))\geq\sum_{s=0}^ka(k,s)^2$, with equality when $r_k$ is irreducible, and non-equivalent to $r_1,\ldots,r_{k-1}$. But $\sum_{s=0}^ka(k,s)^2=\dim(Mor(v^k,v^k))$, and by Corollary 2, we have equality.

Finally, since any irreducible representation of $A_o(F)$ must appear in some tensor power of $u$, and we have a formula for decomposing each $u^k$ into sums of representations $r_s$, we conclude that these representations $r_s$ are all the irreducible representations of $A_o(F)$.
\end{proof}

\medskip

{\bf Remarks.} (1) By recurrence on $k$, we have $\dim(r_k)=(x^{k+1}-y^{k+1})/(x-y)$, where $x,y$ are the solutions of $X^2-nX+1=0$. At $n=2$ we have $\dim(r_k)=k+1$.

(2) The proof of Theorem 1 shows that the commutant of $u^k$ is precisely $A_{\beta,k}$. This can be used in order to prove some simplicity results, in the spirit of \cite{hsk}.

\medskip

Finally, Theorem 1 has the following converse, which follows from the definition of $A_o(F)$, and from Theorem 1 itself:

\begin{theorem}
If the irreducible representations of a compact quantum group $G$ are self-adjoint, and can be indexed by $\mathbb N$, with $r_0=1$, $r_1=u$ and
$r_kr_s=r_{|k-s|}+r_{|k-s|+2}+\ldots+r_{k+s-2}+r_{k+s}$, then $G_{red}$ must be similar to a certain $A_o(F)_{red}$.
\end{theorem}

{\bf Acknowledgements.} I would like to thank G. Skandalis, my PhD advisor, and E. Blanchard, for numerous useful discussions on the subject.

\end{document}